\documentclass{amsart}
\usepackage[utf8]{inputenc}

\usepackage{amsmath}
\usepackage{graphicx}
\usepackage{scalerel}
\usepackage{stackengine,wasysym}
\usepackage{hyperref}

\usepackage{amssymb}
\usepackage{amsthm}
\usepackage{mathtools}
\usepackage{adjustbox}
\usepackage{dsfont}
\usepackage{MnSymbol}
\usepackage{mathrsfs}
\usepackage{stmaryrd}
\usepackage[all]{xy}
\usepackage[mathcal]{eucal}
\usepackage{verbatim}  
\usepackage{fullpage}  
\usepackage{hyperref}
\usepackage{setspace}
\usepackage{tikz-cd}
\usepackage{xcolor}
\usetikzlibrary{positioning}
\renewcommand{\phi}{\varphi}

\usepackage{rotating}
\usepackage{amssymb}

\theoremstyle{definition}

\numberwithin{equation}{subsection}
\begin{document}
\title{Remarks on Hilbert's Thirteenth Problem\\
(Bemerkungen Zum Dreizehnten Hilbertschen Problem)}
\author{Ludwig Bieberbach\\
Translator: Anubhav Nanavaty}
\date{Received January 5, 1931}
\maketitle
\begin{abstract}
    This is an English translation of Ludwig Bieberbach's paper ``Remarks on Hilbert's Thirteenth Problem" originally written in German and originally published in Journal für die Reine und Angewandte Mathematik - 165 (89-92) 1931, along with an addendum to the paper published in 1933. Beiberbach studies under what conditions are there functions of three variables which can or cannot be obtained by combining or nesting functions of two variables. In the addendum, Bieberbach acknowledges a fatal error in a proof in the original article and connects Hilbert's thirteenth problem to a related problem in the theory of differential equations.
\end{abstract}
In 1900, Hilbert posed 23 problems\begin{footnote}{The lecture is not only included in the Comptes rendus of the Paris International Mathematical Congress,
but has also been reprinted in Göttinger Nachrichten 1900 and in Archiv für öfathematik und Physik.}\end{footnote}. They have almost all been solved in the last
30 years\begin{footnote}{In``Die Naturewissenschaften" 1930, Book 61, I gave an overview of this.}\end{footnote}.

The thirteenth of these problems, however, has generally been passed by carelessly. The purpose of the following note is to show that this is not due to the inaccessibility of these questions, and that with relatively easy effort one can gain something from this
problem.

The problem refers to the representation of functions of three variables by nesting
functions of two variables. If $x_1,x_2,x_3$ are three variables, then $f(x_1,x_2),g(x_2,x_3),h(x_3,x_1)$
are functions of two variables. If $a,b,c,A,B,C,\dots$ are other functions of two variables, e.g. $a(f, g), b(g, h), c(h,f), A (a, b), B(b, c), C(A,B),...$ are functions of three variables. Functions which can be represented by using finitely
many functions of two variables are obtained, we will say, by nesting functions of two variables. In his lecture, Hilbert sets the task of proving that a certain function he names
cannot be obtained in this way. He further remarks that he has convinced himself by
rigorous consideration that there are analytic functions of three variables which cannot be obtained ``by finitely-multiple concatenation of functions of only two arguments". Hilbert's statement must first be correctly understood to mean that the used functions of two arguments are also analytic. Hilbert indeed occasionally proved in lectures that there are entire functions of three variables which cannot be represented by nesting analytic functions of two variables. If one wanted to take Hilbert's assertion as literally as it is stated in his lecture, it would not even be correct. In fact, Hilbert himself has remarked in lectures that any function of three variables can be obtained by interpolating functions of two variables, if one means the function concept in Dirichlet's sense. Hereafter I will
show that there are continuous functions of three variables which cannot be obtained
by telescoping continuous functions of two variables. So, the restrictions, which one
imposes on the function term, are the ones which give content to the problem.
\section{Each function of three variables can be expressed by nesting functions of two variables}

For example, let $f(x_1,x_2,x_3)$ be defined in the cube $0\leq x_i\leq 1$. Suppose $0\leq x_2\leq 1$, $0\leq x_3\leq 1$ maps bijectively to the distance $0\leq y\leq 1$  and set $x_1=x$. Then the cube is bijectively mapped to the square $0\leq x\leq 1, 0\leq y\leq 1$. Let $x_1=x$, $x_2=\alpha(y)$, $x_3=\beta(y)$ be the mapping. Then, 
\[f(x_1,x_2,x_3)=f(x,\alpha(y),\beta(y))=F(x,y)\]
Further, the mapping can also be represented as:
\[x=x_1,\quad y=\phi(x_2,x_3)\]
So,
\[f(x_1,x_2,x_3)=F(x_1,\phi(x_2,x_3))\]
Also it follows that every function of three variables is representable by two functions of two variables. However, discontinuous functions certainly occur.
\section{Not every continuous function of three variables can be represented by a combination of continuous functions of three variables}
The proof is based on the fact that any continuous function of two variables can be
uniformly approximated by polynomials of two variables.\\ I call a function \emph{representable} if it can be represented by nesting finitely many continuous functions of two variables.\\
I call a function \emph{$n$-representable}, if it can be represented by $n$ continuous functions of two variables at most. So every representable function is also $n$-representable for large enough $n$.\\ \emph{$n$-rational} or \emph{$n$-polynomial} is a polynomial which can be represented by at most $n$ polynomials of two variables.\\
A function is called \emph{$n$-approximable} if it can be uniformly approximated by $n$-rational functions in a given domain.\\
Every $n$-representable function is $n$-approximable. This follows from Weieratraß's theorem, according to which any function of two variables which is continuous in a closed domain can be approximated uniformly by polynomials of two variables. This is probably clear without further ado if the functions of two variables we have used are continuous at all values of these two variables. By the theorem to be proved, however, it is only required that the functions in question are continuous on certain closed sets of points which are defined by the set of values of the
functions to be used. As the range of values of the three variables $x_1,x_2,x_3$ whose continuous functions are concerned, we take a closed area $B$ such that $x_1,x_2,x_3\in B$. The Tietze extension theorem \begin{footnote}{H. Tietse, Über Functionen, die auf einer abgeschlossenen öfenge atetig eind. Crelles Journal vol.
146 (1916), p. M23. Cf. also F. Hausdorff, Über halbstetige Funktionen und deren Verallgemeinerung. Math. Ztschr. Bd. 5 (1919) Pg 232-239. - Addendum: Just now, Herr stud. math. Rado, in Berlin, has just published a proof which establishes an explicit representation of the extension function in few lines.}\end{footnote} teaches, however, that every function that is continuous on a closed set can be extended to a function that is continuous in the whole plane and that agrees with the given one on the closed set.\\
\emph{Each non-n-approximable function has a neighborhood of non-n-approximable functions.}\\
If $f(x_1,x_2,x_3)$ is not $n$-approximable, Let $f_i$ be a sequence of $n$-approximable functions. In $B$ let:
\[|f_i-f|<\frac{\epsilon_i}{2},\quad \epsilon_i>0,\quad \epsilon_i\to 0,\quad i\to \infty\]
Let $g_i$ denote a sequence of $n$-rational functions and let $|f_i-g_i|<\frac{\epsilon}{2}$ in $B$. Then, $|f-g_i|<\frac{\epsilon}{2}$ in $B$. If all this were true, then $f$ would be n-approximable. Therefore there is an $\epsilon > 0$ such that no $F$ is n-approximable if $|f-F|<\epsilon$.\\
\emph{There are non-$k$-polynomials for each $k$}. If $f(x_1,x_2,x_3)$ is rational, then also $f(x_1,x_2,x_3)-f(0,0,0)$ is rational. So let $f(0,0,0)=0$ and $f(x_1,x_2,x_3)=\phi(f_1,f_2)$ where $f_1$ and $f_2$ are $k-1$ rational. So,
\begin{align*}
    \phi(f_1,f_2)&=\phi\{f_1(0,0,0)+f_1'(x_1,x_2,x_3),f_2(0,0,0)+f_2'(x_1,x_2,x_3)\}\\
    &=\Phi(f_1',f_2'),\quad\quad \Phi(0,0)=0
\end{align*}
Here $f_1',f_2'$ are new $k-1$ polynomials with $f_1'(x_1,x_2,x_3)=f_2'(x_1,x_2,x_3)=0$. So we may assume that all polynomials used for the representation of $k$-polynomials $f(x_1,x_2,x_3)$ vanish at the
origin. Then to get those $f(x_1,x_2,x_3)$ that grow at most to the $r$-th order in each of the variables, one only needs the terms up to the $r$-th order in each variable for the polynomials of two variables used in the representation. At most $k$ such polynomials are used. Each yields at most $r^3$ coefficients. $f(x_1,x_2,x_3)$ has $r^3$ coefficients in its given form. But since $r^3>kr^2$ for $r>k$, it is not possible to specify the $k$-many polynomials of two variables in such a way that determines the corresponding polynomial of three variables.\\
\emph{Thus, for every degree exceeding $k$, there are non-$k$ polynomials. The coefficients of the $k$-polynomials whose degree exceeds $k$ must satisfy certain algebraic conditional equations. Among these equations are some that are independent of the degree}.\\
\emph{For all $k$ there are polynomials of three variables which are not $k$-approximizable}.\\ Let us choose a $k+1$-degree polynomial which is not $k$-polynomial, i.e. whose coefficients do not obey one of the degree-independent conditional equations. Let the polynomial be $f(x_1,x_2,x_3)$. Then there is an $\epsilon>0$, so that the coefficients
of any polynomial $F$, for which $|f-F|<\epsilon$, satisfy the respective conditional equation. This is because, from $|f- F| < \epsilon$, estimates for the differences of $f$ and $F$ in each coefficient follow, and these estimates approach zero with $\epsilon$. (The conditional equations refer to the coefficients of the members of order at most $k+1$). So there is a neighborhood of $f$ that does not contain a $k$-polynomial. Thus, $f$ is also not $k$-approximable.\\
\emph{In the neighborhood of every continuous function $g(x_1,x_2,x_3)$ there are, for every $n$, continuous functions that are not $n$-approximable}.\\
If $g$ is not $n$-approximable, then the assertion is trivial. However, if $g$ itself is $n$-approximable, let $f(x_1,x_2,x_3)$ be a function that is not $2n$-approximable. Then for every $t$, $tf$ is also not $2n$-approximable. Then
\[h=g+tf\]
is not $n$-approximable, because otherwise $tf-h-g$ would be $2n$-approximable. One can choose $f$ as a polynomial. Thus, if g is a polynomial, then h is also a polynomial.\\
\emph{There are non-representable continuous functions of three variables}\begin{footnote}
    {Translator's Note: This is false; see the Kolmogorov--Arnold Theorem}
\end{footnote}. $f_1(x_1,x_2,x_3)$ is not $n_1$-approximable, and $\epsilon_1>0$ is chosen such that for any $F$, such that $|f_1-F|<\epsilon_1$ in $B$, is not $n_1$ approximable.  Let $f_2$ be chosen from this neighborhood such that it is $n_2$-approximable for $n_2>n_1$. Let $\epsilon_2$ be chosen such that no $F$ is $n_2$-approximable for $|f_2-F|< \epsilon$, and so on for each $i$. Let $n_i\to\infty$ and $\epsilon_i\to 0$. Then the $f_i(x_1,x_2,x_3)$ in B converge uniformly to a limit function $f(x_1,x_2,x_3)$
which is not $n$-approximizable for any $n$. Therefore, this function is also no longer $n$-representable for any n. So it is not
representable at all. In the preceding proof all $f_i$ can be chosen as polynomials.
\section{Not every analytical function of three real variables can be represented by a combination of continuous functions of two variables\begin{footnote}
    {Translator's Note: This is false; see the Kolmogorov--Arnold Theorem}
\end{footnote}}
To see this, one has only to take care that the polynomials constructed in the last part of the just finished proof converge uniformly in a certain range of the three complex variables. For this purpose, the bounds $|f_i-F|<\epsilon_i$ used there are always related to this fixed range, to which the otherwise used range $B$ of the real variables belongs as a part. Then the $f_i$ are uniformly convergent in this range. The limit function is then not representable. But it is also analytic at the same time. By slight modification of the construction one can even achieve that the limit function is an entire function.
\section{}
The method given here can be applied to many other problems. Ostrowski has shown in Math. Ztschr. vol. 8 (1920) that the function
\[\zeta(u,v)=\sum_{n=1,2\dots}\frac{u^n}{n^v}\]
cannot be represented by using finitely many analytic functions of one variable and
arbitrary algebraic processes in finite number. Hilbert recently (Math. Ann. 97 (1927)) referred to these questions again and in particular posed the problem to investigate to what extent functions of two variables can be represented by using finitely many functions of one variable and by using the sum process finitely often. Also here the answer depends on the used notion of function. If one thinks of analytic functions, then the mentioned result of Ostrowski contains the answer. If one thinks of continuous functions, ao the method of this work allows to give the answer. If one finally takes the notion of function in the Dirichletian sense\begin{footnote}
    {In the sense of Dirichlet, these functions are set-functions, i.e. there is no continuity assumption.}
\end{footnote}, then \emph{it can be shown that every function of two variables can be represented by finitely many functions of one variable and finitely often by using the summation process.} Perhaps there will be another opportunity to go into this not quite trivial remark in more detail.\\
Received January 5, 1931.
\end{document}


\title{Addendum to my Work ``Remarks on Hilbert's Thirteenth Problem" in Volume 165 of This journal}
\author{Ludwig Bieberbach\\
Translator: Anubhav Nanavaty}
\date{29 October 1933}
\maketitle

On p. 91 there is an error. The proof that for each $k$ there are $k$-polynomials of three variables which are not $k$-approximable is not valid. This error in my proof, which was more based on hastiness than ignorance, was brought to my attention by Mr. Kamke. My efforts to repair the damage have not yet led to the desired result. In the meantime, I have found an approach that has some prospects and that I would like to share here because it is related to another question that is also interesting to me and apparently and strangely has not been dealt with so far.

This is the question of whether there exist differential equations (in our case partial
algebraic differential equations) by whose solutions arbitrary continuous functions can
be approximated uniformly (in the function values, not also in the derivatives).

The following consideration leads to the connection of the Hilbert problem with this
question:

In development of the considerations of the second paragraph on p. 91 of my work one
shows that every $r\geq 3k-3$ times differentiable $k$-function (in particular therefore every $k$-polynomial) satisfies an $r$-th order algebraic partial differential equation. If then it can be shown that not every continuous function of three variables can be uniformly approximated by solutions of such a differential equation, then it was shown that there are $k$-polynomials of
three variables which are not $k$-approximable. With this, the gap of the proof would be closed.